\documentclass[12pt]{article} 
\usepackage[hidelinks]{hyperref}
\usepackage[all,cmtip]{xy}
\usepackage{color}
\usepackage{stackrel}
\usepackage{amsmath,amssymb,eucal,bbm} 
\font\tenrm=cmr10

\font\cmssl=cmss10 at 12 pt  
   
\font\bigss=cmssdc10 scaled 2300

\font\cmsslll=cmss10 at 14 pt

\renewcommand{\a}{\alpha}  
\renewcommand{\b}{\beta}  
  
\renewcommand{\d}{\delta}  
\newcommand{\e}{\epsilon}  
\newcommand{\f}{\varphi}  
\newcommand{\g}{\gamma}

\renewcommand{\o}{\omega}  
  
\renewcommand{\r}{\rho}  
\newcommand{\s}{\sigma}  
\renewcommand{\t}{\tau}

\newcommand{\z}{\zeta}  
  
\newcommand{\F}{\Phi}  
\newcommand{\G}{\Gamma}

\renewcommand{\O}{\Omega}

\newcommand{\bC}{\mathbb{C}}  
\newcommand{\bR}{\mathbb{R}}  
\newcommand{\bZ}{\mathbb{Z}}  
\newcommand{\bH}{\mathbb{H}}

\newcommand\SO{\mathrm{SO}}

\newcommand{\p}{\partial}  
    
\renewcommand{\square}{\kern1pt\vbox  
               {\hrule height 0.6pt\hbox{\vrule width 0.6pt\hskip 3pt  
    \vbox{\vskip 6pt}\hskip 3pt\vrule width 0.6pt}\hrule height0.6pt}  
    \kern1pt}  

\newcommand{\ra}{\rightarrow}

\DeclareMathOperator\End{End\;}

\newtheorem{Pb}{Problem}

\newtheorem{Th}{Theorem}  
\newtheorem{Prop}{Proposition}  
  
\newtheorem{Cor}{Corollary}  
\newtheorem{Lem}{Lemma}  
\newtheorem{Def}{Definition} 
\newcommand{\bP}{\begin{Pb}\ \ } 
\newcommand{\eP}{\end{Pb}}  
\newcommand{\bt}{\begin{Th}\ \ }  
\newcommand{\et}{\end{Th}}  
\newcommand{\bp}{\begin{Prop}\ \ }  
\newcommand{\ep}{\end{Prop}}  
\newcommand{\bc}{\begin{Cor}\ \ }  
\newcommand{\ec}{\end{Cor}}  
\newcommand{\bl}{\begin{Lem}\ \ }  
\newcommand{\el}{\end{Lem}}  
\newcommand{\bd}{\begin{Def}\ \ }  
\newcommand{\ed}{\end{Def}}  
  
\newcommand{\pf}{\noindent{\it Proof:\ \ }}  
\newcommand{\qed}{\hfill\square}  
\newcommand{\n}{\nabla}  
  
\newcommand{\ot}{\otimes}

\newcommand{\be}{\begin{equation}}  
\newcommand{\ee}{\end{equation}}  
  
\newcommand\re[1]{(\ref{#1})}  
\newcommand{\arr}{\begin{array}{rlll}}  
\newcommand{\ea}{\end{array}}  
\newcommand{\bea}{\begin{eqnarray}}  
\newcommand{\eea}{\end{eqnarray}}  
\newcommand{\bean}{\begin{eqnarray*}}  
\newcommand{\eean}{\end{eqnarray*}}  
  
\catcode`@=11  
\@addtoreset{equation}{section}  
\catcode`@=12  
 
 \newcommand{\zb}{\bar{z}}
\newcommand{\wb}{\bar{w}}
\renewcommand{\z}{\zeta}
\newcommand{\zt}{\tilde{\zeta}}
\renewcommand{\F}{\mathtt{F}}
\newtheorem{Rem}{Remark}  
\newcommand{\br}{\begin{Rem}\ \ }  
\newcommand{\er}{\end{Rem}}

\begin{document}  
 \rightline{} 
\vskip 1.5 true cm  
\begin{center}  
{\bigss  Quaternionic K\"ahler metrics associated with\\[1em] 
special K\"ahler manifolds}\\[.5em]
\vskip 1.0 true cm   
{\cmsslll  D.V.\ Alekseevsky$^1$, V.\ Cort\'es$^{2,3}$, M.\ Dyckmanns$^2$ and T.\ 
Mohaupt$^4$} \\[3pt] 
$^1${\tenrm Institute for Information Transmission Problems\\
B.\ Karetny per.\ 19\\ 
127051 Moscow, Russia\\
and\\ 
Masaryk University\\
Kotlarska 2\\ 
61137 Brno, Czech Republic\\
dalekseevsky@iitp.ru}\\[1em]  
{\tenrm   $^2$Department of Mathematics\\  
University of Hamburg\\ 
Bundesstra{\ss}e 55, 
D-20146 Hamburg, Germany\\
cortes@math.uni-hamburg.de\\ 
malte.dyckmanns@math.uni-hamburg.de}\\[1em]  
{\tenrm   $^3$Center for Mathematical Physics\\ 
University of Hamburg\\ 
Bundesstra{\ss}e 55, 
D-20146 Hamburg, Germany }\\[1em]
$^4${\tenrm Department of Mathematical Sciences\\ 
University of Liverpool\\
Peach Street, Liverpool L69 7ZL, UK\\  
Thomas.Mohaupt@liv.ac.uk}\\[1em]   
March 4, 2015
\end{center}  
\vskip 1.0 true cm  
\baselineskip=18pt  
\begin{abstract}  
\noindent  
We give an explicit formula for the quaternionic K\"ahler metrics 
obtained by the HK/QK correspondence.  As an application, 
we give a new proof of the fact that the Ferrara-Sabharwal metric
as well as its one-loop deformation is quaternionic K\"ahler.
A similar explicit formula is given for the analogous
(K/K) correspondence between K\"ahler manifolds endowed with a
Hamiltonian Killing vector field. As an example, we apply this formula 
in the case of an arbitrary conical K\"ahler manifold.    
\end{abstract}

\tableofcontents
\section*{Introduction}
Extending results by Haydys \cite{Haydys}, it was proven in \cite{ACM} that any pseudo-hyper-K\"ahler manifold  $(M,g,J_1,J_2,J_3)$ 
of dimension $4n$ endowed with a 
space-like or time-like $\o_1$-Hamiltonian Killing vector field $Z$ which acts as a rotation in the plane 
spanned by  $J_2$ and $J_3$ 
gives rise to a one-parameter family of conical\footnote{See
  Definition \ref{conicalDef}.} pseudo-hyper-K\"ahler manifolds of dimension 
$4n+4$ and finally 
to a one-parameter family of possibly indefinite quaternionic K\"ahler manifolds of dimension $4n$. 
Here $\o_\a := g J_\a := g\circ J_\a := g(J_\a \cdot ,\cdot )$,  $\a = 1,2,3$, are the three symplectic forms 
associated with the pseudo-hyper-K\"ahler structure and  
the parameter in the above one-parameter families is related to the choice of a Hamiltonian function for $Z$.  
Under the assumptions on the Hamiltonian specified in \cite{ACM}, the resulting quaternionic K\"ahler metrics are positive definite. 

Following \cite{APP,Hi}  (but allowing indefinite metrics) we will call the above relation between  hyper-K\"ahler and quaternionic K\"ahler manifolds
of the same dimension the HK/QK correspondence. The analogous construction
relating (possibly indefinite) K\"ahler  manifolds of the same dimension, which follows from the K\"ahler conification in 
\cite{ACM},  will be called the K/K correspondence.  

It was also proven in \cite{ACM}  that 
the cotangent bundle of any conical special K\"ahler manifold admits a 
canonical vector field $Z$ 
which satisfies the above assumptions with respect to the pseudo-hyper-K\"ahler 
structure $(g,J_1,J_2,J_3)$ provided by the (rigid) c-map \cite{CFG} (see Section \ref{sectionRigidCmap}). 
Using techniques 
from supergravity and twistor theory, Alexandrov, Persson and Pioline 
\cite{APP} show that  
the Ferrara-Sabharwal metric \cite{FS} (also known as the 
supergravity c-map metric, see Section \ref{sugraSect}) and its
one-loop deformation are related to the c-map pseudo-hyper-K\"ahler
metric $g$ under the HK/QK correspondence. 
It was shown in \cite{ACM} that the above vector field $Z$ has a canonical  
Hamiltonian function conjecturing 
 that the quaternionic K\"ahler metric associated with 
this particular choice of the parameter is precisely the 
Ferrara-Sabharwal metric. It was checked
that the sign of the scalar curvature is negative and thus consistent
 with the latter conjecture. 
Finally, the precise relation between the parameter in the choice of the 
Hamiltonian and the 
one-loop quantum deformation 
parameter occurring in \cite{RSV,APP} was left for future investigation.


In this paper we verify the above conjecture and determine  the precise 
relation between the Hamiltonian parameter and the one-loop
parameter. In fact, we apply the HK/QK correspondence to 
the pseudo-hyper-K\"ahler manifolds obtained from the 
rigid c-map starting with a conical affine special
K\"ahler manifold. The final result is the  
formula \eqref{DefFSmetric} for the quaternionic K\"ahler metric, 
see Theorem \ref{HKQKThm}. This is precisely the
one-loop deformed Ferrara-Sabharwal metric as described in 
\cite{RSV,APP}. 
As a corollary this implies: 
\bc The Ferrara-Sabharwal metric and its
one-loop deformation \eqref{DefFSmetric} are quaternionic K\"ahler.
\ec 
Notice that this generalizes the result that
the Ferrara-Sabharwal metric is quaternionic K\"ahler \cite{FS,Hi1}. 

Our proof is based on a new explicit formula for
the quaternionic K\"ahler metric in the  HK/QK correspondence, see
Theorem \ref{mainThm}. A similar result is obtained in the K\"ahler case, that is for the 
K/K correspondence, see
Theorem \ref{mainThmKK}. 
To obtain the explicit formula for the quaternionic K\"ahler metric
we start by reviewing the Swann bundle construction and the moment map
of a tri-holomorphic Killing vector field on the Swann bundle in Section 
\ref{SwannSec}. 
Our approach allows to control the signature
of the resulting metrics. In particular, we specify
for any given value of the one-loop parameter $c$  
the maximal domain on which the deformed Ferrara-Sabharwal
metric is positive definite. For $c\ge 0$, this domain coincides with
the manifold on which the Ferrara-Sabharwal
metric is defined. These results generalize 
those of Antoniadis, Minasian, Theisen and Vanhove \cite{AMTV} 
in four dimensions (for the universal
hypermultiplet). 

We have included 
appendix A, in which we discuss the simplest case of the 
HK/QK correspondence in which the 
initial hyper-K\"ahler manifold is (flat) four-dimensional, 
for the reader's convenience. The resulting 
quaternionic K\"ahler manifold is the complex
hyperbolic plane (universal hypermultiplet).  

For the K/K correspondence we apply our formula
for the resulting metric in the case when the initial
pseudo-K\"ahler manifold is conical, see Theorem \ref{thmKKexample}. In particular,
for a conical affine special K\"ahler manifold
$(M,J,g,\n ,\xi )$ we obtain (up to a cyclic covering)
the product $\bC H^1\times \bar{M}$ of the projective
special K\"ahler manifold  $\bar{M}$ underlying $M$ 
and the complex hyperbolic line, see Remark \ref{RemarkKK}. Notice that this 
is a maximal totally geodesic K\"ahler submanifold
$\bC H^1\times \bar{M}\subset \bar{N}$ 
of the Ferrara-Sabharwal manifold $\bar{N}$, which is related to  
$\bar{M}$ by the supergravity c-map.

 \noindent
{\bf Acknowledgments} 
This work was part of a research 
project within the RTG 1670 ``Mathematics inspired by String Theory'',
funded by the Deutsche Forschungsgemeinschaft (DFG). D.V.A. has been
supported by the project CZ.1.07/2.3.00/20.0003 of the Operational
Programme Education for Com\-petitiveness of the Ministry of
Education, Youth and Sports of the Czech Republic. The work of T.M. was
supported in part by STFC grant ST/G00062X/1. 
We thank Nigel Hitchin, 
Andrew Swann and Owen Vaughan for discussions.

\section{The Swann bundle revisited}\label{SwannSec}
In this section we derive explicit formulas relating the 
metric of a quaternionic 
K\"ahler manifold to the pseudo-hyper-K\"ahler metric of its Swann bundle \cite{Swann}. 
This will be used in Section \ref{ExplHKQKSec} to obtain an explicit formula for the 
quaternionic K\"ahler metric in the HK/QK correspondence
from the conical pseudo-hyper-K\"ahler metric constructed in \cite{ACM}. 
\subsection{The pseudo-hyper-K\"ahler structure}\label{SwannbdlSec}
Let $(M,g,Q)$ be a (possibly indefinite) quaternionic K\"ahler manifold of
nonzero scalar curvature, where $Q\subset \mathfrak{so}(TM)$ denotes its quaternionic structure. 
Let us denote by $\pi : S\ra M$ the principal $\mathrm{SO}(3)$-bundle
of frames $(J_1,J_2,J_3)$ in $Q$ such that $J_3=J_1J_2$ and $J_\a^2=-\mathrm{Id}$, $\a =1,2,3$.   
The principal action of  an element $A\in \mathrm{SO}(3)$ is given by 
\[ s= (J_1,J_2,J_3)\mapsto \t (A,s) := (J_1,J_2,J_3)A^\e,\]
where $\e=1$ if we consider $S$ as a right-principal bundle and $\e=-1$ if we
prefer a left-principal bundle. 
Let us denote by $Z_\a$ the fundamental vector fields associated with 
some basis $(e_\a)$ of $\mathfrak{so}(3)$: 
\[ Z_\a (s) =\left. \frac{\partial}{\partial t}\right|_{t=0}\t (\exp (te_\a ),s).\]  
We may choose 
the basis corresponding to the standard basis  of
$\mathfrak{sp}(1)=\mathrm{Im} \bH\cong \bR^3$ under the canonical 
isomorphism $\mathfrak{sp}(1) \cong ad (\mathfrak{sp}(1))=\mathfrak{so}(3)$. 
Then 
\be \label{so3Equ} [e_\a ,e_\b ] = 2e_\g, \quad [Z_\a ,Z_\b ] = 2\e Z_\g,  \ee
for every cyclic permutation $(\a ,\b ,\g )$ of $(1,2,3)$. 
In the following, $(\a ,\b ,\g )$ will be always a cyclic permutation, whenever
the three letters appear in an expression. 

The Levi-Civita connection $\n$ of $(M,g)$ induces a principal connection
\[ \theta =\sum \theta_\a e_\a : TS \ra \mathfrak{so}(3)\] 
on $S$. Its curvature is defined by 
\[ \O := d\theta +\e \frac{1}{2}[\theta \wedge \theta] ,\]
where 
\[ \frac{1}{2}[\theta \wedge \theta] (X,Y) := [\theta(X) , \theta (Y)],\quad X,Y\in T_sS,\quad s\in S.\]
Writing $\O = \sum \O_\a e_\a$ and using \re{so3Equ} we have 
\be \label{OaEqu} \O_\a = d\theta_\a +2\e \theta_\b\wedge
\theta_\g.\ee 
From the definition of the connection and curvature
forms we get the following lemma.
\bl \label{LthetaLemma} 
\[ \mathcal{L}_{Z_\a}\theta_\a = \mathcal{L}_{Z_\a}\O_\a =0,\quad
\mathcal{L}_{Z_\a}\theta_\b = 2\e \theta_\g,\quad 
\mathcal{L}_{Z_\a}\O_\b = 2\e \O_\g .\]

\el

Given a local section $\boldsymbol{\s} = (J_1,J_2,J_3) \in \G (U,S)$, defined over some open subset $U\subset M$, 
we can also define a vector-valued 1-form 
\[ \bar{\theta}=\sum \bar{\theta}_\a e_\a\]
on $U$ by 
\[ \n J_\a = -2\e (\bar{\theta}_\b \ot J_\g -\bar{\theta}_\g \ot J_\b) .\]
The coefficient is chosen such that 
\[ \n (J_1,J_2,J_3) = (J_1,J_2,J_3)\e  \bar\theta .\]
Notice that then 
\be \label{covEqu} \n B = dB+ \e  \sum \bar{\theta}_\a \ot [ J_\a, B],\ee
for every section $B= \sum b_\a J_\a$ of $Q$, where $d=d_{\boldsymbol{\s}}$ is defined by $dB :=\sum db_\a \ot J_\a$. The vector-valued 1-forms $\bar{\theta} $ on $U\subset M$ and
$\theta$ on $S$ are related by 
\[ \bar{\theta} = \boldsymbol{\s}^*\theta .\]
In the local trivialization $\pi^{-1}(U)\cong U\times \mathrm{SO}(3)$ of $S$
given by $\s$ we can write 
\[ \quad \theta = \pi^*\bar{\theta} + \varphi,\]
where $\f  = \sum \f_\a e_\a$ is the Maurer-Cartan form on $ \mathrm{SO}(3)$ defined
by $\f_\a (Z_\b) = \d_{\a \b}$.  
{}From \re{covEqu} we compute the curvature $R^Q\in \G (\wedge^2T^*M\ot Q)$,
$Q \cong ad (Q)\subset \End Q$,  of the vector bundle $Q$, which 
is   
\[ R^Q = \sum \bar{\O}_\a J_\a,\quad \bar{\O}_\a = \e d\bar{\theta}_\a +2\bar{\theta}_\b\wedge \bar{\theta}_\g.\] 
It is a well-known result by Alekseevsky \cite{A}  that  
\[ \bar\O_\a = -\frac{\nu}{2} \o_\a,\]
where $\o_\a = gJ_\a$ and 
\[ \nu := \frac{scal}{4n(n+2)}\quad  (\dim M=4n)\]
is the reduced scalar curvature. Since the curvature form
of a principal connection is horizontal, this implies that
\be \label{OEqu} \e \O_\a\big|_{\boldsymbol{\s} (U) = \pi^* \bar{\O}_\a\big|_{\boldsymbol{\s} (U)} = -\frac{\nu}{2} \pi^*\o_\a\big|_{\boldsymbol{\s} (U)}.}\ee

\noindent 
We endow the manifold $S$ with the pseudo-Riemannian metric
\[ g_S = \sum \theta^2_\a + \frac{\nu}{4} \pi^*g.\] 

\noindent
Now we consider the cone $\hat{M}= S\times \bR^{>0}$ over $S$ with 
the Euler vector field $\xi := Z_0 := r\p_r$ and the 
following exact 2-forms 
\[ \hat{\o}_\a :=  d\hat{\theta}_\a,\quad \hat{\theta}_\a:=-\e
\frac{r^2}{2}\theta_\a .\]
For later use we state the following lemma, which follows from 
Lemma \ref{LthetaLemma} and the fact that $Z_0=\xi$ preserves $\theta_\a$.
\bl The Lie algebra 
$\mathrm{span}\{ Z_i| i= 0,\ldots,3\}\cong \mathfrak{co}(3)$ acts on
$\mathrm{span}\{ \hat{\theta}_\a|\a =1,2,3\}$ by the 
standard representation:
\[  \mathcal{L}_{Z_0}\hat{\theta}_\a = 2\hat{\theta}_\a,\quad
\mathcal{L}_{Z_\a}\hat{\theta}_\b = 2\e \hat{\theta}_\g .\]
\el 
Using the above data we recover Swann's hyper-K\"ahler structure on $\hat{M}$: 

\bt The cone metric $\hat{g} = dr^2 + r^2g_S$  is 
a pseudo-hyper-K\"ahler metric on $\hat{M}$ with the K\"ahler forms  $\hat{\o}_\a$.  
The signature of $\hat{g}$ is $(4+4k,4l)$ if $\nu >0$ and $(4+4l,4k)$ if $\nu <0$,
where $(4k,4l)$ is the signature of the quaternionic K\"ahler metric $g$ on $M$. 
\et 

\pf Let us denote by $T^v\hat{M}\subset T\hat{M}$ the vertical distribution with respect to
the projection $\hat{\pi}:= \pi \circ \mathrm{pr}_S :\hat{M} \ra M$, $\mathrm{pr}_S : \hat{M} = S\times \bR^{>0}\ra S$,
and by $T^h\hat{M}$ the horizontal 
distribution defined by its $\hat{g}$-orthogonal complement. 
Let $\hat{J}_\a$ be the uniquely determined 3 almost complex structures on $\hat{M}$  which 
preserve the horizontal distribution and satisfy   
\[ \hat{J}_\a Z_0= -\e Z_\a,\quad \hat{J}_\a Z_\a = \e Z_0,\quad  \hat{J}_\a Z_\b = Z_\g,\quad \hat{J}_\a Z_\g = -Z_\b, 
\quad  \hat{\pi}_* \circ \hat{J}_\a|_{(s,r)} =J_\a \circ \hat{\pi}_*, \]
where $s=(J_1,J_2,J_3)$. We see that
these structures satisfy $\hat{J}_1\hat{J}_2=\hat{J}_3$ and pairwise anti-commute. 
Then,  using \re{OaEqu} and \re{OEqu}, one can easily check 
$\hat{g}\hat{J}_\a = \hat{\o}_\a$. 
This proves that the  
2-forms $\hat{\o}_\a$ are not only closed but also non-degenerate
and that $\hat{J}_\a=-\hat{\o}_\b^{-1}\hat{\o}_\g$ are three anti-commuting skew-symmetric almost  
complex structures on $(\hat{M},\hat{g})$.  
By the Hitchin Lemma \cite[Lemma 6.8]{Hi0}, this shows that $(\hat{g},\hat{J}_1,\hat{J}_2,\hat{J}_3)$ is a pseudo-hyper-K\"ahler
structure on $\hat{M}$. 
\qed 
\subsection{The moment map of an infinitesimal automorphism}
\label{momentSection}
Let $\hat{M}$ be the Swann bundle over a (possibly indefinite) quaternionic K\"ahler manifold $(M,g,Q)$. 
We will follow the conventions in Section \ref{SwannbdlSec} with $\e=-1$. We endow $\hat{M}$  with the hyper-K\"ahler structure 
$(g_{\hat{M}} := \s \hat{g},(\hat{J}_\a ) )$, where $\s = \pm 1$. The corresponding 
K\"ahler forms are $\s \hat{\o}_\a=d(\s \hat{\theta}_\a)$. 

Let $X$ be a tri-holomorphic 
space-like or time-like Killing vector field on $\hat{M}$, which
commutes 
with the 
Euler vector field $\xi = r\p_r=Z_0$. 
\bp The vector field $X$ is tri-Hamiltonian with moment
map $-\mu$, where
\[ \mu :\hat{M} \ra \bR^3,\quad 
x\mapsto (\mu_1(x),\mu_2(x),\mu_3(x)),\quad \mu_\a := \hat{\theta}_\a
(X).\] 
 In fact, the functions $\mu_\a$ satisfy
\be \label{muEqu} d\mu_\a = -\iota_X\hat{\o}_\a .\ee
\ep 

\pf Notice first that, since $X$ is tri-holomorphic, it commutes  not only with $\xi$ but also with $Z_\a = \hat{J}_\a \xi$.
This implies already that the Killing field $X$ preserves the horizontal distribution $T^h\hat{M}=(T^v\hat{M})^\perp$
and hence the three one-forms $\theta_\a$. Furthermore,  $\mathcal{L}_X(r^2)=\mathcal{L}_Xg_{\hat{M}}(\xi,\xi)=0$, since $X$ is Killing and commutes with $\xi$. This implies that 
\[ \mathcal{L}_X\hat{\theta}_\a=\mathcal{L}_X\left( \frac{r^2}{2}\theta_\a \right)=0.\] 
Using this equation, we have
\[ d\mu_\a = d\iota_X\hat{\theta}_\a =  \mathcal{L}_X\hat{\theta}_\a - \iota_Xd\hat{\theta}_\a = 
-\iota_X\hat{\o}_\a.\]
\qed 

We will now explain how to recover the 
quaternionic K\"ahler metric on $M$
from the geometric data on the level set of the moment map
$\mu$ 
\[ P=\{ \mu_1 = 1,\; \mu_2=\mu_3=0\}\subset \hat{M}.\] 
Since the group $\bR^{>0}\times \SO(3)$ generated by 
$\xi, Z_1, Z_2,Z_3$ acts as the standard conformal linear group $\mathrm{CO}(3)$ on the three-dimensional vector space spanned by the 
functions $\mu_\a$, i.e.\
\[ \mathcal{L}_{Z_0}\mu_\a = 2\mu_\a,\quad
\mathcal{L}_{Z_\a}\mu_\b = -2 \mu_\g ,\]
we  see that 
\[ \hat{M}\setminus \{ \mu =0\}  = \bigcup_{a\in \bR^{>0}\times \SO(3)}aP . \] 
In particular, $P$ is nonempty. Then \re{muEqu} shows that 
$P\subset \hat{M}$ is a smooth submanifold of codimension 3. 
On $P$ we have the following data:
\begin{eqnarray*}
g_P &:=&  g_{\hat{M}}|_P=\s \hat{g}|_{P}\in \G (\mathrm{Sym}^2T^*P)\\
\theta_\a^P & := & \s \hat{\theta}_\a|_P\in \O^1(P)\;  (\a=1,2,3)\\
f &:=& \s \left. \frac{r^2}{2}\right|_P \in C^\infty (P)\\
\theta_0^P  &:=& -\frac12 df\in \O^1(P)\\
X_P &:=& \sigma X|_P \in \mathfrak{X}(P)\\
Z_1^P &:=&  Z_1|_P\in \mathfrak{X}(P). 
\end{eqnarray*} 
The fact that $Z_1$ is tangent to $P$ follows from 
\[ d\mu_\a Z_1 = \iota_{Z_1}d\mu_\a = \mathcal{L}_{Z_1}\mu_\a = -2\d_{2\a }\mu_3 +2\d_{3\a} \mu_2,\]
since $\mu_2=\mu_3=0$ on $P$.

With these definitions, the formula 
\be \label{ghatEqu}\hat{g} = dr^2 + r^2\left(\sum_{\a =1}^3 \theta_\a^2 + \frac{\nu}{4}\pi^*g\right)\ee
implies:
\bp \label{auxProp} The quaternionic K\"ahler metric $g$ on $M$ is related as follows to the 
geometric data on the level set $P\subset \hat{M}$ of the moment map: 
\be  \label{PropEqu}\nu\pi^*g|_P = \frac{2}{f}\left( g_P - \frac{2}{f}\sum_{a=0}^3 (\theta_a^P)^2\right). \ee

\ep

\pf Solving Eq.\ \re{ghatEqu} for  $\nu\pi^*g$ yields:
\[ \nu \pi^*g= \frac{4}{r^2}(\hat{g} -dr^2 -r^2\sum\theta_\a^2)=\frac{4}{\s r^2}(\s \hat{g} -\s dr^2 -\s r^2\sum\theta_\a^2).\] 
Restricting to $P$ we first obtain:
\be \label{firstEqu}\nu \pi^*g|_P= \frac{2}{f}(g_P -\s dr^2|_P -2f\sum\theta_\a^2|_P).\ee
The above definitions imply $\theta_\a|_P =  f^{-1}\theta_\a^P$. Therefore, 
$f\theta_\a^2|_P =f^{-1}(\theta_\a^P)^2$. 
Similarly, $\s dr^2|_P= 2f^{-1}(\theta_0^P)^2$. This shows that \re{firstEqu} implies \re{PropEqu}. 

\qed 

\bc \label{auxCor} The tensor field on the right-hand side of  \re{PropEqu} is invariant under $Z_1^P$ and has one-dimensional kernel
$\bR Z_1^P$. 
\ec 
\pf The $Z_1^P$-invariance follows from the $Z_1$-invariance of
$\pi^*g$. The statement about the kernel follows from
\be \label{intersecEqu}T^v\hat{M}\cap TP = \bR Z_1,\ee 
which is a consequence of
\begin{eqnarray*} d\mu_\a \xi  &=& \mathcal{L}_\xi \mu_\a = 2 \mu_\a \\
d\mu_\a Z_2 & =& \mathcal{L}_{Z_2}\mu_\a = 2\d_{1\a} \mu_3 -2\d_{3\a} \mu_1\\
d\mu_\a Z_3 &=& \mathcal{L}_{Z_3}\mu_\a = -2\d_{1\a} \mu_2 +2\d_{2\a} \mu_1.
\end{eqnarray*}
In fact, we have already shown that the vertical vector field $Z_1$ is tangent to $P$ and these equations show
now that the  three vector fields $\xi, Z_2, Z_3$  are 
mapped to (constant) linearly independent vectors under the vector-valued one-form $d\mu = (d\mu_\a)  : T\hat{M}\ra \bR^3$.  
This implies \re{intersecEqu}, since $TP = \ker d\mu$. 
\qed 

In the next section we apply the above results to the case when $\hat{M}$ 
is obtained by conification of a hyper-K\"ahler manifold, in the sense of \cite{ACM}.
 \section{Explicit formula for the HK/QK correspondence}
\label{ExplHKQKSec}
Let $(M,g,J_1,J_2,J_3)$ be a possibly indefinite hyper-K\"ahler manifold
with the K\"ahler forms $\o_\alpha = gJ_\alpha$, $\alpha =1,2,3$, and a 
time-like or space-like $\o_1$-Hamiltonian Killing vector field $Z$ such that 
$\mathcal{L}_ZJ_2=-2J_3$. 
According to \cite{ACM}, 
with any choice of function $f\in C^\infty (M)$ such that 
$df = -\o_1 Z$ and such that $f_1 = f-\frac{g(Z,Z)}{2}$ is not 
zero, one can, at least locally, associate 
a quaternionic K\"ahler metric $g'$ on a manifold
$M'$ of dimension $\dim M$. (One has to assume, 
in particular,  that
the functions $f$ and $f_1$ are nowhere zero, which may require 
to restrict the manifold $M$.)

Following \cite{ACM}, let  $P\ra M$ be an $S^1$-principal bundle 
with a principal connection $\eta$ with the curvature
$d\eta = \o_1 -\frac{1}{2}dgZ$. We endow $P$ with the pseudo-Riemannian metric 
\be \label{gPEqu}g_P := \frac{2}{f_1}\eta^2 + \pi^*g\ee 
and with the vector field 
\be \label{Z1PEqu} Z_1^P :=\tilde{Z} +f_1X_P,\ee 
where $\tilde{Z}$ denotes the horizontal lift of $Z$ and
$X_P$ the fundamental vector field of the principal action.    
Furthermore, we endow $P$ with the following 
one-forms: 
\begin{eqnarray}
\theta_0^P  &:=& -\frac{1}{2}df\nonumber \\ 
\theta_1^P &:=& \eta +\frac{1}{2}gZ\nonumber \\ 
\theta_2^P  &:=& \frac{1}{2}\o_3Z\nonumber \\ 
\theta_3^P  &:=& -\frac{1}{2}\o_2Z.\label{thetaEqu} 
\end{eqnarray}

\noindent 
\bt \label{mainThm} The tensor field 
\be \label{gPtildeEqu} \tilde{g}_P:=g_P-\frac{2}{f}\sum_{a=0}^3(\theta_a^P)^2\ee
on $P$ is invariant under $Z_1^P$ and has one-dimensional kernel $\bR Z_1^P$. Let $M'$ be a codimension $1$ submanifold of $P$ which 
is transversal to
the vector field $Z_1^P$. Then 
\[ g' := \frac{1}{2|f|}\tilde{g}_P|_{M'}\] 
is a possibly indefinite
quaternionic K\"ahler metric on $M'$.  
\et 

\pf Analysing the 
constructions of \cite{Haydys,ACM}, we see that the data 
$g_P$, $\theta_\a^P$ , $f$, $X_P$, $Z_1^P$ are obtained 
by restriction from data $\s \hat{g}$, $\s \hat{\theta}_\a$, $\s r^2/2$, $X$ and $Z_1$ on the conical pseudo-hyper-K\"ahler manifold
$\hat{M}$, as in Section \ref{momentSection}. Therefore, the theorem follows from Proposition \ref{auxProp} and Corollary \ref{auxCor} . The tensor field $\frac{1}{2|f|}\tilde{g}_P$
corresponds to $\frac{\s \nu}{4} \pi^*\bar{g}|_P $, where $(\bar{M},\bar{g})$ denotes the underlying 
quaternionic K\"ahler manifold, when $\hat{M}$ is represented locally as a Swann bundle
$\hat{\pi} : \hat{M} \ra \bar{M}$. (Recall that $\s = \mathrm{sgn} f$.) 
\qed 

The  metric $g'$ is the quaternionic K\"ahler metric 
which corresponds under the 
HK/QK correspondence to the hyper-K\"ahler manifold $(M,g,(J_\a))$ 
with the data $(Z,f)$. Notice that the principal projection $\pi : (S,g_{\hat{M}}\big|_S=\sigma g_S)\ra (\bar{M},\bar{g})$ is a pseudo-Riemannian
submersion if and only if $\frac{\sigma \nu }{4}=1$.
This is why we normalized the metric $g'$ such that
its reduced scalar curvature is $\nu'=4\sigma$.

\noindent 
\br \rm{If $Z_1^P$ generates a free and proper 
action of a one-dimensional Lie group
 $A$ ($\cong S^1$ or $\bR$) and if $M'$ is a global section for the $A$-action, then 
 we can identify $M'$ with the orbit space $P/A$, which inherits the  
 quaternionic K\"ahler metric $g'$.} 
\er 
 
 In the next section we present a similar result for the K/K correspondence.  
 \section{Explicit formula for the K/K correspondence}
 Let $(M,g,J)$ be a possibly indefinite K\"ahler manifold endowed with a 
time-like or space-like  Killing vector field $Z$, which is Hamiltonian
with respect to the  K\"ahler form $\o = gJ$. According to \cite{ACM}, 
with any choice of function $f\in C^\infty (M)$ such that 
$df = -\o Z$ and such that $f_1 = f-\frac{g(Z,Z)}{2}$ is not 
zero, one can, at least locally, associate a conical 
pseudo-K\"ahler manifold  $\hat{M}$ of (real) dimension $\dim M +2$ 
and, hence,  a 
pseudo-K\"ahler manifold $M'$ of dimension $\dim M$. 
In fact, $\hat{M}$ is a metric cone
over a pseudo-Sasaki manifold $S$ which has a 
pseudo-K\"ahler structure transversal to the Reeb foliation. 
Therefore, any codimension $1$ submanifold of $S$  transversal to the Reeb foliation 
inherits a pseudo-K\"ahler structure $(J',g')$. Now we give an
explicit formula for the metric $g'$ in terms of the initial data. 

Following \cite{ACM}, let  $\pi: P\ra M$ be an $S^1$-principal bundle 
with a principal connection $\eta$ with the curvature
$d\eta = \o -\frac{1}{2}dgZ$. We endow $P$ with the pseudo-Riemannian metric 
\[  g_P := \frac{2}{f_1}\eta^2 + \pi^*g\]
and with the vector field 
\[ Z^P :=\tilde{Z} +f_1X_P,\]
where $\tilde{Z}$ denotes the horizontal lift of $Z$ and
$X_P$ the fundamental vector field of the principal action.    
Furthermore, we endow $P$ with the following 
one-forms: 
\begin{eqnarray*}
\theta_0^P  &:=& -\frac{1}{2}df\nonumber \\ 
\theta_1^P &:=& \eta +\frac{1}{2}gZ .\nonumber 
\end{eqnarray*}
Then $\hat{M}=\bR \times P$ is endowed 
with a conical pseudo-K\"ahler structure 
described explicitly in terms of the above data on $P$, see \cite{ACM}. 
In particular, the Euler vector field is given by $\xi = \p_t$, where $t$ is the 
coordinate on the $\bR$-factor.  It is related to the radial variable $r>0$ 
of the metric cone over the pseudo-Sasaki manifold $S$ by $e^{2t}=\frac{r^2}{2|f|}$.  
This implies that $S= \{ p\in \hat{M}| r(p) =1\}$  is a circle bundle over $M$ diffeomorphic to $P$. 

\noindent 
\bt \label{mainThmKK} The tensor field 
\[ \tilde{g}_P:=g_P-\frac{2}{f}\sum_{a=0}^1(\theta_a^P)^2\]
on $P$ is invariant under $Z^P$ and has one-dimensional kernel $\bR Z^P$. Let $M'$ be a codimension $1$ submanifold of $P$ which 
is transversal to
the vector field $Z^P$. Then 
\[ g' := \frac{1}{2|f|}\tilde{g}_P|_{M'}\] 
is a possibly indefinite 
K\"ahler metric on $M'$.  
\et 

\pf The proof is similar to that of  Theorem \ref{mainThm}. 
It relies on the representation of the pseudo-K\"ahler manifold $(\hat{M}, g_{\hat{M}})$ as a 
metric cone over a pseudo-Sasaki manifold $S$.  $\hat{M} = \bR^{>0}\times S$ is equipped 
with the metric $g_{\hat{M}}= \s \hat{g}=\s (dr^2 + r^2g_S)$, where $\s = \mathrm{sgn} f\in \{ -1, 1\}$. 
One can (locally) assume that $S= I\times \bar{M}\subset \bR \times \bar{M}$ is contained in  a trivial principal bundle
with structure group $\bR$ over a pseudo-K\"ahler manifold $(\bar{M},\bar{g})$, where $I\subset \bR$ is an interval. 
Let us denote by $\bar{\o}$ the K\"ahler form of $(\bar{M},\bar{g})$. 
The pseudo-Sasaki metric takes the form $g_S = \theta^2 + \bar{g}$, where $\theta$ is a principal connection with curvature
given by $2\bar{\omega}$.
Analysing the 
construction of \cite{ACM}, we see that the tensor field $\frac{1}{2|f|}\tilde{g}_P$
corresponds to $\s \pi^*\bar{g}|_P $, where $\pi : \hat{M} \ra \bar{M}$ is the composition of the two projections $\hat{M}\ra S$ and 
$S\ra \bar{M}$. Here $P=\{ t=0\}\times P\subset \hat{M}= \bR \times P$ is the level set $\{ \mu =1\}$ 
of the moment map $\mu = e^{2t}$ associated with the holomorphic Killing  vector field $X$ on $\hat{M}$ which 
canonically extends the vector field $X_P$ on $P$. 
\qed 
\renewcommand{\d}{\partial}
\subsection{K/K correspondence for conical K\"ahler manifolds}
As an example, we apply the K/K correspondence to an arbitrary
conical pseudo-K\"ahler manifold $(M,J,g,\xi)$ endowed with the
holomorphic Killing field $Z:=2J\xi$. Recall the following
definition:
\bd \label{conicalDef} A pseudo-Riemannian
manifold $(M,g)$ is called {\cmssl conical} if it is endowed 
with a space-like or time-like vector field $\xi$ (called the {\cmssl
  Euler
vector field}) such that
$D\xi = \mathrm{Id}$, where $D$ is the Levi-Civita connection. 
\ed 
Geometrically this means that $M$ is locally isometric to 
a (space-like or time-like, 
respectively) metric cone $C_\pm (S)= (\bR^{>0}\times S,
\pm dr^2+r^2g_S)$ over a pseudo-Riemannian manifold $(S,g_S)$. 
Notice that in this local representation the Euler vector field $\xi$
is given by $r\p_r$. If $g$ happens to be a pseudo-K\"ahler metric
for some complex structure $J$ on $M$, then $(M,J,g,\xi )$ is called
a  {\cmssl conical pseudo-K\"ahler manifold}. In this case  $M$ is
locally isometric to a {\cmssl pseudo-K\"ahler cone}, that
is a metric cone $C_\pm (S)$ over a pseudo-Sasaki manifold $(S,g_S)$ 
with Reeb vector field $J\xi|_S$, see e.g.\ \cite{BC,MSY}. 

From now on
we assume that $(M,J,g,\xi)$ is a conical pseudo-K\"ahler manifold. 
Using $r^2:=|g(\xi,\xi)|$, $\lambda:=\mathrm{sgn}\,g(\xi,\xi)$, $\tilde{\eta}:=\frac{\lambda}{r^2}g(J\xi,\cdot)$, we can write the metric as \begin{equation} g=\frac{(g(\xi,\cdot))^2}{g(\xi,\xi)}+\frac{(g(J\xi,\cdot))^2}{g(\xi,\xi)}+|g(\xi,\xi)|\,\breve{g}=\lambda (dr^2+r^2(\tilde{\eta}^2+\lambda \breve{g})).\label{KconeMetric}\end{equation}
This equation defines the tensor $\breve{g}$ on $M$, which has $\ker \breve{g}=\mathrm{span}\{\xi,J\xi\}$ and fulfills $\mathcal{L}_\xi\breve{g}=\mathcal{L}_{J\xi}\breve{g}=0$. Assume that $S:=\{r=1\}\subset M$ is non-empty and let $\breve{M}\subset S$ be a codimension $1$ submanifold that is transversal to the Reeb vector field $J\xi|_S\in \Gamma (TS)$. Then $\breve{M}$ inherits a complex structure $\breve{J}$ from $J$ such that $(\breve{M},\breve{J},\breve{g}|_{\breve{M}})$ is pseudo-K\"ahler.
For simplicity (and without restriction of generality), we assume in the following theorem
that $M=\mathbb{R}^{>0}\times S$ is globally a cone. 

\bt\label{thmKKexample}
The K/K correspondence assigns to any pseudo-K\"ahler cone \linebreak ($M=\mathbb{R}^{>0}\times S,J,g,\xi$) endowed with the holomorphic Killing field $Z=2J\xi$ the manifolds $M'_\pm:=I_\pm\times S^1\times\breve{M}$,
\[I_+:=\begin{cases} (\max\{0,-2c\},\infty) &~\mbox{for } \lambda=1 \\
(\min\{-2c,0\},-c) &~\mbox{for } \lambda=-1,\end{cases}\quad 
I_-:=\begin{cases} (-c,\max\{0,-2c\}) &~\mbox{for } \lambda=1 \\
(-\infty,\min\{-2c,0\}) &~\mbox{for } \lambda=-1;\end{cases}\]
 endowed with the metric
\begin{equation}
g'=\frac{1}{2|\rho|}\left[\lambda(\rho+c)\breve{g}-\frac{1}{4\rho}\frac{\rho+2c}{\rho+c}d\rho^2-\frac{1}{4\rho}\frac{\rho+c}{\rho+2c}(d\tilde{\phi}-2c\tilde{\eta}|_{\breve{M}})^2\right] \label{KKExampleMetric}
\end{equation}
for each $c\in \mathbb{R}$. Here, $\tilde{\phi}$ is a local coordinate on the $S^1$-factor $S^1=\{e^{-\frac{i}{4}\tilde{\phi}}|\tilde{\phi}\in\mathbb{R}\}$, $\rho\in I_\pm$ and $\lambda=\mathrm{sgn}\,g(\xi,\xi)$. The signature of $(M'_+,g')$ is $(2k,2l+2)$ and that of $(M'_-,g')$ is $(2k+2,2l)$, where $(2k,2l)$ is the signature of $\breve{g}$.\\
For $c=0$, we get \[(M'_\pm,2g')=(\pm \mathbb{R}^{>0}\times S^1\times
\breve{M},\mp g_{\mathbb{C}H^1}+\breve{g})\cong (\mathbb{R}^{>0}\times S^1\times
\breve{M},\mp g_{\mathbb{C}H^1}+\breve{g}),\] 
where $\pm$ corresponds to $\lambda=\pm 1$, respectively, and $g_{\mathbb{C}H^1}:=\frac{1}{4\rho^2}(d\rho^2+d\tilde{\phi}^2)$.
\et
\pf
$f:=\lambda r^2-c$ fulfills $\omega(Z,\cdot)=-2g(\xi,\cdot)=-2\lambda r \,dr=-df$, where $r^2=|g(\xi,\xi)|$. The K\"ahler form of $(M,J,g)$ is given by $\omega=\lambda r dr\wedge \tilde{\eta}+r^2\breve{\omega}$, where $\breve{\omega}:=\breve{g}(J\cdot,\cdot)$. One can check that $\breve{\omega}=\frac{1}{2}\lambda d\tilde{\eta}$. Using this, one finds $d\beta=4\omega$, where $\beta:=g(Z,\cdot)=2\lambda r^2\tilde{\eta}$.

We endow the trivial $S^1$-principal bundle $P:=M\times S^1\to M$ with the principal connection \[\eta:=ds-\frac{1}{4}\beta=ds-\frac{\lambda}{2}r^2\tilde{\eta},\] which has curvature $d\eta=\omega-\frac{1}{2}d\beta=-\omega$. Here, $s$ is the natural coordinate on $S^1=\{e^{is}|s\in\mathbb{R}\}$. The metric and one-forms on $P$ are given by
\begin{align*}
g_P&=\frac{2}{f_1}\eta^2+g\\
\theta_0^P&=-\frac{1}{2}df=-\lambda r\,dr\\ \theta_1^P&=\eta+\frac{1}{2}\beta=ds+\frac{\lambda}{2}r^2\tilde{\eta},
\end{align*}
where $f_1=f-\frac{1}{2}g(Z,Z)=-\lambda r^2-c$.

We compute the degenerate tensor field $\tilde{g}_P$:
\begin{align*}
\tilde{g}_P=g_P-\frac{2}{f}((\theta^P_0)^2+(\theta_1^P)^2)&=\frac{2}{f_1}(ds+\frac{f_1+c}{2}\tilde{\eta})^2+g-\frac{2}{f}(r^2\,dr^2+(ds+\frac{f+c}{2}\tilde{\eta})^2)\\
&\!\!\stackrel{\re{KconeMetric}}{=}\left(\frac{2}{f_1}-\frac{2}{f}\right)(ds+\frac{c}{2}\tilde{\eta})^2+(\lambda-\frac{2r^2}{f})dr^2+r^2\breve{g}\\
&=-\frac{4}{f}\frac{f+c}{f+2c}(ds+\frac{c}{2}\tilde{\eta})^2-\frac{1}{4f}\frac{f+2c}{f+c}df^2+\lambda(f+c)\breve{g}.
\end{align*}  

Since $\mathbb{R}^{>0}\times \breve{M} \subset \mathbb{R}^{>0} \times S= M$ is transversal to $J\xi\in \Gamma(TM)$, $\tilde{M}':=\mathbb{R}^{>0}\times \breve{M}\times S^1\subset P$ is transversal to $Z^P:=\tilde{Z}+f_1 \partial_s=Z-(\eta(Z)-f_1)\partial_s=2J\xi-c\partial_s\in \Gamma(TP)$. Replacing the coordinates $r$ and $s$ by $\rho:=f$ and $\tilde{\phi}:=-4s$, we obtain the K\"ahler metric $g'=\frac{1}{2|\rho|}\tilde{g}_P|_{M'}$ obtained from the K/K correspondence (Theorem \ref{mainThmKK}) as given in Eq.\ \eqref{KKExampleMetric}. Here,
\[M':=\begin{cases} (-c,\infty) \times S^1\times\breve{M} &\quad\quad\mbox{for } \lambda=1 \\
(-\infty,-c) \times S^1\times\breve{M} &\quad\quad\mbox{for } \lambda=-1
\end{cases}
\]
is obtained from $\tilde{M}'$ via the coordinate change $r\mapsto \rho=\lambda r^2-c$. For the metric $g'$ to be defined, we need to restrict to $\{f=\rho\neq 0,~ -f_1=\rho+2c\neq 0\}\subset M'$.

The signature of $g$ is given by $(2k+2,2l)$ if $\lambda=1$ and $(2k,2l+2)$ if $\lambda=-1$, where $(2k,2l)$ is the signature of $\breve{g}$. The signature of $g'$ is  related to the one of $g$ by
\[\mathrm{sign}\,g'=\begin{cases} (+2,-2)+\mathrm{sign}\,g &\quad\quad\mbox{for } f_1>0,~f<0 \\
\quad\quad\quad\quad\quad \mathrm{sign}\,g &\quad\quad\mbox{for } ff_1>0
\\
(-2,+2)+\mathrm{sign}\,g &\quad\quad\mbox{for } f_1<0, f>0.
\end{cases}\]
Using $f=\rho$, $f_1=-(\rho+2c)$ and taking into account $r^2=\lambda(\rho+c)>0$, one finds that on the subsets $M_\pm=\{\rho\in I_\pm\}\subset M'$ given in the Theorem, $g'$ has signature $(2k,2l+2)$, $(2k+2,2l)$ respectively.

For the last statement, one just has to notice that for $c=0$, $\mathrm{sgn}\,\rho=\lambda$.
\qed
\section{HK/QK correspondence for the c-map}
In this section, we use the explicit formula given in Theorem
\ref{mainThm} to show that the pseudo-hyper-K\"ahler structure on the
cotangent bundle of a conical affine special K\"ahler manifold given
by the rigid c-map is related to the quaternionic K\"ahler metric
obtained from the supergravity c-map via the HK/QK correspondence. In
fact, we get a one-parameter family of positive definite quaternionic
K\"ahler metrics, which corresponds to one-loop corrections of the
hypermultiplet moduli space in string theory compactifications on
Calabi-Yau 3-folds (if the corresponding model is realized in string
theory). 
As a corollary, this proves that the Ferrara-Sabharwal metric and its one-loop deformation are indeed quaternionic K\"ahler.

\subsection{Conical affine and projective special K\"ahler geometry}
\label{ConicalSect}
First, we recall the definitions of conical affine and projective special K\"ahler manifolds \cite{ACD,CM}:
\bd
A {\cmssl conical affine special K\"ahler manifold} $(M,J,g_M,\nabla,\xi)$ is a pseudo-K\"ahler manifold $(M,J,g_M)$ endowed with a flat torsionfree connection $\nabla$ and a vector field $\xi$ such that
\begin{enumerate}
\item[i)] $\nabla \omega_M=0$, where $\omega_M:=g_M(J\cdot,\cdot)$ is the K\"ahler form,
\item[ii)] $(\nabla_XJ)Y=(\nabla_YJ)X$ for all $X,Y\in \Gamma(TM)$,
\item[iii)] $\nabla\xi=D\xi=\mathrm{Id}$, where $D$ is the Levi-Civita connection,
\item[iv)] $g_M$ is positive definite on $\mathcal{D}=\mathrm{span}\{\xi,J\xi\}$ and negative definite on $\mathcal{D}^\perp$.
\end{enumerate}
\ed

Let $(M,J,g_M,\nabla,\xi)$ be a conical affine special K\"ahler manifold of complex dimension $n+1$. Then $\xi$ and $J\xi$ are commuting holomorphic vector fields that are homothetic and Killing respectively \cite{CM}. We assume that the holomorphic Killing vector field $J\xi$ induces a free $S^1$-action and that the holomorphic homothety $\xi$ induces a free $\mathbb{R}^{>0}$-action on $M$. Then $(M,g_M)$ is a metric cone over $(S,g_S)$, where $S:=\{p\in M|g_M(\xi(p),\xi(p))=1\}$, $g_S:=g_M|_S$; and $-g_S$ induces a Riemannian metric $g_{\bar{M}}$ on $\bar{M}:=S/S^1_{J\xi}$. $(\bar{M},-g_{\bar{M}})$ is obtained from $(M,J,g)$ via a 
K\"ahler reduction with respect to $J\xi$ and, hence, $g_{\bar{M}}$ is a K\"ahler metric (see e.g. \cite{CHM}). The corresponding K\"ahler form $\o_{\bar{M}}$ is
obtained from $\o_M$ by symplectic reduction. This determines
the complex structure $J_{\bar{M}}$. 

\bd
The K\"ahler manifold $(\bar{M},J_{\bar{M}},g_{\bar{M}})$ is called a {\cmssl projective special K\"ahler manifold}.
\ed

More precisely, $S$ is a (Lorentzian) Sasakian manifold and introducing the radial coordinate $r:=\sqrt{g(\xi,\xi)}$, we can write the metric on $M$ as  
\cite{BC,MSY}
\begin{equation}\label{dKconeMetricDecomposition}
g_M=dr^2+r^2 \pi^\ast g_S, \quad g_S=g_M|_S= \tilde{\eta}\otimes \tilde{\eta}|_S-
\bar{\pi}^\ast g_{\bar{M}},
\end{equation}
where
\begin{equation}\label{contactOneForm}
\tilde{\eta}:=\frac{1}{r^2}g_M(J\xi,\cdot)=d^c \log r=i(\overline{\d}-\d)\log r
\end{equation}
is the contact one-form form when restricted to $S$ and $\pi :M\to S=M/\mathbb{R}^{>0}_\xi$, $\bar{\pi}: S\to \bar{M}=S/S^1_{J\xi}$ are the canonical projection maps. From now on, we will drop $\pi^\ast$ and $\bar{\pi}^\ast$ and identify, e.g., $g_{\bar{M}}$ with a $(0,2)$ tensor field on $M$ that has the 
distribution  
$\mathcal{D}=\mathrm{span}\{\xi,J\xi\}$ as its kernel.

Locally, there exist so-called \emph{conical special holomorphic coordinates} $z=(z^I)=(z^0,\ldots ,z^n):U\stackrel{\sim}{\to} \tilde{U}\subset \mathbb{C}^{n+1}$ such that the geometric data on the domain $U\subset M$ is encoded in a holomorphic function $F: \tilde{U}\to \mathbb{C}$ that is homogeneous of degree 2 \cite{ACD,CM}. Namely, we have \cite{CM} \[g_M|_U=\sum_{I,J} N_{IJ}dz^Id\zb^J, \quad N_{IJ}(z,\zb):=2\mathrm{Im}\,F_{IJ}(z):=2\mathrm{Im}\,\frac{\d^2F(z)}{\d z^I \d z^J} \quad (I,J=0,\ldots,n)\]
and $\xi|_U=\sum z^I\frac{\d}{\d z^I}+\zb^I\frac{\d}{\d\zb^I}$. The K\"ahler potential for $g_M|_U$ is given by $r^2|_U=g_M(\xi,\xi)|_U=\sum z^IN_{IJ} \zb^J$.

The $\mathbb{C}^\ast$-invariant functions $X^\mu:=\frac{z^\mu}{z^0}$, $\mu=1,\ldots, n$, define a local holomorphic coordinate system on $\bar{M}$. The K\"ahler potential for $g_{\bar{M}}$ is $\mathcal{K}:=-\log \sum_{I,J=0}^n X^IN_{IJ}(X)
\bar{X}^J$, where $X:=(X^0,\ldots, X^n)$ with $X^0:=1$.

\subsection{The rigid c-map}\label{sectionRigidCmap}
Now, we introduce the \emph{rigid c-map}, which assigns to each affine special (pseudo-)K\"ahler manifold $(M,J,g_M,\nabla)$ and in particular to any conical affine special K\"ahler manifold $(M,J,g_M,\nabla,\xi)$ of real dimension $2n+2$ a (pseudo-)hyper-K\"ahler manifold $(N=T^*M,g_N$, $J_1,J_2,J_3)$ of dimension $4n+4$ \cite{CFG,ACD}.

From now on, we assume for simplicity that $(M\subset \mathbb{C}^{n+1},J=J_{can},g_M,\nabla,\xi)$ is a conical affine special K\"ahler manifold that is globally described by a homogeneous holomorphic function $F$ of degree $2$ 
defined on a $\bC^*$-invariant domain $M$ in standard holomorphic coordinates $z=(z^I)=(z^0,\ldots,z^n)$ induced from $\mathbb{C}^{n+1}$. Here, $J_{can}$ denotes the standard complex structure induced from $\mathbb{C}^{n+1}$.

The real coordinates $(q^a)_{a=1,\ldots,2n+2}:=(x^I,y_J)_{I,J=0,\ldots, n}:=(\mathrm{Re}\,z^I,\,\mathrm{Re}\,F_J(z):=\mathrm{Re}\,\frac{\d F(z)}{\d z^J})$ on $M$ are $\nabla$-affine and fulfill $\omega_M=-2\sum dx^I\wedge dy_I$, where $\omega_M=g(J\cdot,\cdot)$ is the K\"ahler form on $M$ \cite{CM}. We consider the cotangent bundle $\pi_N: N:=T^\ast M\to M$ and introduce real functions $(p_a):=(\zt_I,\z^J)$ on $N$ such that together with $(\pi_N^\ast q^a)$, they form a system of canonical coordinates.

\bp 
In the above coordinates $(z^I,p_a)$, the hyper-K\"ahler 
structure on $N=T^\ast M$ obtained from the rigid c-map is given by
\begin{align}
g_N&=\sum dz^IN_{IJ}d\zb^J+\sum A_IN^{IJ}\bar{A}_J,\label{RigidCmapMetric}\\
\omega_1&=\frac{i}{2}\sum N_{IJ}dz^I\wedge d\zb^J+\frac{i}{2}\sum N^{IJ}A_{I}\wedge \bar{A}_J,\label{RigidCmapOmega1}\\
\omega_2&=-\frac{i}{2}\sum (d\zb^I\wedge \bar{A}_I-dz^I\wedge A_I),\label{RigidCmapOmega2}\\
\omega_3&=\frac{1}{2}\sum (dz^I\wedge A_I+d\zb^I\wedge \bar{A}_I) \label{RigidCmapOmega3},
\end{align}
where $A_I:=d\zt_I+\sum_J F_{IJ}(z)d\zeta^J ~(I=0,\ldots, n)$ are
complex-valued one-forms on $N$ and 
$\omega_\alpha=g_N(J_\alpha\cdot,\cdot)$. 
(Here and in the following, we
identify functions and one-forms on $M$ with their pullbacks to $N$.)
\ep 
\pf 
One can check by a direct calculation that the metric and K\"ahler forms,
\eqref{RigidCmapMetric}--\eqref{RigidCmapOmega3} agree with the
geometric data for the rigid c-map given in Section 3 of
\cite{ACD} (see also Section 3 of \cite{ACM}), up to a conventional sign in the definition of 
the K\"ahler forms $\omega_\alpha=g_N(J_\alpha\cdot,\cdot) = -
g_N(\cdot , J_\a\cdot )$ in \cite{ACD}.
For instance, we can write $\omega_1$ and $\omega_3$ as
\begin{align}
\omega_1&=-2\sum dx^I\wedge dy_I+\frac{1}{2}\sum d\zt_I\wedge d\z^I,\label{RigidCmapOmega1v2}\\
\omega_3&=\sum dx^I\wedge d\zt_I+\sum dy_I\wedge d\z^I=\sum dq^a\wedge dp_a. \nonumber
\end{align}
\qed

\br
\rm{It follows from the intrinsic geometric description in \cite{ACD} that 
the pseudo-hyper-K\"ahler structure is independent of the
particular description of the special K\"ahler structure
in terms of a holomorphic function $F$.} 
\er
\br
\rm{We introduce holomorphic functions $w_I$, $I=0,\ldots,n$, on $(N,J_1)$ that together with the holomorphic coordinates $z=(z^I)$ on $(M,J)$ form a system of canonical holomorphic coordinates on $(N=T^\ast M, J_1)$. Then $(w_I)$ and $(\zt_I,\z^J)$ are related by
\begin{align*} \sum_I w_Idz^I+\wb_Id\zb^I&\stackrel{!}{=} \sum_I \zt_Idx^I+\z^I dy_I\\&=\sum_I \frac{\zt_I}{2}(dz^I+d\zb^I)+\frac{\zeta^I}{2}(\sum_J F_{IJ}(z)dz^J+\overline{F_{IJ}(z)}d\zb^J),\end{align*}
which is equivalent to
\begin{equation} w_I=\frac{1}{2}(\zt_I+\sum_J F_{IJ}(z)\z^J)\quad\quad (I=0,\ldots, n). \label{RelFiberCoordsCmap}\end{equation}
With the identification \eqref{RelFiberCoordsCmap}, \eqref{RigidCmapMetric}--\eqref{RigidCmapOmega3} also agree, up to conventional factors, with the rigid c-map as given in Appendix B of \cite{CFG} and throughout the physics literature.}
\er

\subsection{The supergravity c-map}
\label{sugraSect}
Let $(\bar{M},g_{\bar{M}})$ be a projective special K\"ahler manifold of complex dimension $n$ which is globally defined by a single holomorphic function $F$. The \emph{supergravity c-map} \cite{FS} associates with $(\bar{M},g_{\bar{M}})$ 
a quaternionic K\"ahler manifold $(\bar{N},g_{\bar{N}})$ of dimension $4n+4$. 
Following the conventions of \cite{CHM}, we have $\bar{N}=\bar{M}\times \mathbb{R}^{>0}\times \mathbb{R}^{2n+3}$ and
\begin{eqnarray*} \label{e:2} g_{\bar{N}} &=& g_{\bar{M}} + g_G,\\\nonumber 
g_G&=& \frac{1}{4\rho^2}d\rho^2 + \frac{1}{4\rho^2}(d\tilde{\phi}
+ \sum (\zeta^Id\tilde{\zeta}_I-\tilde{\zeta}_Id\zeta^I) )^2 
+\frac{1}{2\rho}\sum \mathcal{I}_{IJ}(m) d\zeta^Id\zeta^J\\ 
&&+ \frac{1}{2\rho}\sum \mathcal{I}^{IJ}(m)(d\tilde{\zeta}_I + 
\mathcal{R}_{IK}(m)d\zeta^K)(d\tilde{\zeta}_J + 
\mathcal{R}_{JL}(m)d\zeta^L),
\end{eqnarray*}
where $(\rho,\tilde{\phi},\tilde{\z}_I, \z^I)$, $I=0,1,\ldots, n$, 
are standard coordinates on $\bR^{>0}\times \bR^{2n+3}$. 
The real-valued matrices $\mathcal{I}(m):=(\mathcal{I}_{IJ}(m))$ and $\mathcal{R}(m):=(\mathcal{R}_{IJ}(m))$
depend only on $m\in \bar{M}$ and $\mathcal{I}(m)$ is invertible with
the inverse $\mathcal{I}^{-1}(m)=:(\mathcal{I}^{IJ}(m))$. More precisely,
 \begin{equation} 
\label{FRI}
{\cal N}_{IJ} := 
\mathcal{R}_{IJ} + i\mathcal{I}_{IJ} := 
\bar{F}_{IJ} + 
i \frac{\sum_K N_{IK}z^K\sum_L N_{JL}z^L}{\sum_{IJ}N_{IJ}z^Iz^J} ,\quad 
N_{IJ} := 2 \mathrm{Im} F_{IJ} ,\end{equation}
where $F$ is the holomorphic prepotential with respect
to some system of special holomorphic coordinates $z^I$ on the 
underlying conical special K\"ahler manifold $M\ra \bar{M}$. 
Notice that the expressions are homogeneous of degree zero and, hence, 
well defined functions on $\bar{M}$. It is shown in \cite[Cor.\ 5]{CHM} 
that the matrix $\mathcal{I}(m)$ is positive definite
and hence invertible and that the metric $g_{\bar{N}}$ does not
depend on the choice of special coordinates \cite[Thm.\ 9]{CHM}. 
It is also shown that $(\bar{N},g_{\bar{N}})$ is complete if 
and only if $(\bar{M},g_{\bar{M}})$ is complete \cite[Thm.\ 5]{CHM}.

Using $(p_a)_{a=1,\ldots, 2n+2}:=(\zt_I,\z^J)_{IJ=0,\ldots,n}$ and $(\hat{H}^{ab}):=\begin{pmatrix}\mathcal{I}^{-1} & \mathcal{I}^{-1}\mathcal{R} \\ \mathcal{R}\mathcal{I}^{-1} & \mathcal{I}+\mathcal{R}\mathcal{I}^{-1}\mathcal{R}\end{pmatrix}$, we can combine the last two terms of $g_G$ into $\frac{1}{2\rho}\sum dp_a \hat{H}^{ab} dp_b$, i.e. the quaternionic K\"ahler metric is given by
\begin{equation}\label{FSmetric} g_{FS}:=g_{\bar{N}}=g_{\bar{M}}+\frac{1}{4\rho^2}d\rho^2 + \frac{1}{4\rho^2}(d\tilde{\phi}
+ \sum (\zeta^Id\tilde{\zeta}_I-\tilde{\zeta}_Id\zeta^I) )^2+ \frac{1}{2\rho}\sum dp_a \hat{H}^{ab} dp_b.\end{equation}

\subsection{HK/QK correspondence for the c-map}
Again, we assume that $(M\subset \mathbb{C}^{n+1},J=J_{can},g_M,\nabla,\xi)$ is a conical affine special K\"ahler manifold that is globally described by a homogeneous holomorphic function $F$ of degree $2$ in standard holomorphic coordinates $z=(z^I)=(z^0,\ldots,z^n)$ induced from $\mathbb{C}^{n+1}$. We want to apply the HK/QK correspondence to the hyper-K\"ahler manifold $(N=T^\ast M,g_N,J_1,J_2,J_3)$ of signature $(4,4n)$ obtained from the rigid c-map (see Section \ref{sectionRigidCmap}). In \cite{ACM}, it was shown that the vector field $Z:=2(J\xi)^h=2J_1\xi^h$ on $N$ fulfills the assumptions of the HK/QK correspondence, i.e. it is a space-like $\omega_1$-Hamiltonian Killing vector field with $\mathcal{L}_{Z}J_2=-2J_3$. Here, $X^h\in \Gamma(TN)$ is defined for any vector field $X\in \Gamma(TM)$ by $X^h(\pi_N^\ast q^a)=\pi_N^\ast X(q^a)$ and $X^h(p_a)=0$ for all $a=1,\ldots, 2n+2$. ($X^h$ is the horizontal lift with respect to the flat
connection $\n$.) 

\bt\label{HKQKThm} 
Applying the HK/QK correspondence to $(N,g_N,J_1,J_2,J_3)$ endowed with 
the $\o_1$-Hamiltonian Killing vector field $Z$ gives (up to a constant
conventional factor) the one-parameter family $g_{FS}^c$ \eqref{DefFSmetric} of 
quaternionic pseudo-K\"ahler metrics, which includes the 
Ferrara-Sabharwal metric $g_{FS}$ \eqref{FSmetric}. The metric
$g_{FS}^c$ is positive definite and of negative scalar curvature
on the domain $\{ \rho > -2c\} \subset \bar{N}$ (which 
coincides with $\bar{N}$ if $c\ge 0$, see Section \ref{sugraSect}). If $c<0$ the metric $g_{FS}^c$
is of signature $(4n,4)$ on the domain $\{ -c < \rho < -2c\} \subset \bar{N}$.
Furthermore, if $c>0$ the metric $g_{FS}^c$
is of signature $(4,4n)$ on the domain $\bar{M} \times \{ -c < \rho < 0\} 
\times \bR^{2n+3}\subset \bar{M} \times \bR^{<0}\times  \bR^{2n+3}$.
\et

\pf
We start from the hyper-K\"ahler structure on $N=T^\ast M$ given in Eqs.\ \eqref{RigidCmapMetric}--\eqref{RigidCmapOmega3}. As in Section \ref{sectionRigidCmap}, we identify functions and differential forms on $M$ with their pullbacks to $\pi_N:N\to M$. We first compute the geometric data involved in the 
HK/QK correspondence, cf.\ Section \ref{ExplHKQKSec}. 
The moment map for $-\omega_1$ w.r.t. $Z=2(J\xi)^h$ is given by $f:=r^2-c$, where $r:=||\xi||_{g_M}=\sqrt{\sum z^I N_{IJ}\zb^J}$ and $c\in\bR$:
\[\omega_1(Z,\cdot)=-g_M(2\xi,\cdot)=-\sum (z^IN_{IJ}d\zb^J+N_{IJ}\zb^J dz^I)=-d(r^2)=-df,\]
since $\sum_I z^I\frac{\partial F_{IJ}(z)}{\partial z^K}=0$. With $g_N(Z,Z)=4g_M(\xi,\xi)=4r^2$, we get
\[f_1:=f-\frac{1}{2}g_N(Z,Z)=-r^2-c.\]
For the functions $f$ and $f_1$ nowhere to vanish, we have to restrict $N$ to $\{r^2\neq |c|\}\subset N$.
Using the contact one form $\tilde{\eta}:=\frac{1}{r^2}g_M(J\xi,\cdot)$ on $M$ (see \eqref{contactOneForm}), we get
\[\beta:=g_N(Z,\cdot)=2g_M(J\xi,\cdot)=2r^2\tilde{\eta}.\]
We consider the trivial $S^1$-principal bundle
\[P:=N\times S^1,\quad S^1=\{e^{is}|s\in \mathbb{R}\},\]
with the connection form
\[\eta=ds+\eta_N,\]
where $\eta_N$ is the following one-form on $N$:
\[\eta_N:=-\frac{1}{2}r^2\tilde{\eta}+\eta_{can}=\frac{f_1+c}{2}\tilde{\eta}+\eta_{can},\quad \eta_{can}:=\frac{1}{4}\sum(\zt_Id\z^I-\z^Id\zt_I).\]
Then
\[d\eta=d\eta_N=-\frac{1}{4}d\beta+d\eta_{can}=\omega_1-\frac{1}{2}d\beta,\]
where we used that $\omega_1$ can be written as \[\omega_1\stackrel{\eqref{RigidCmapOmega1v2}}{=}\pi_N^\ast \omega_M+\frac{1}{2}\sum d\zt_I\wedge d\z^I=\frac{1}{4}d\beta+d\eta_{can},\]
since $\pi_N^\ast \omega_M=\frac{1}{4}\pi_N^\ast dd^c(r^2)$ and $\pi_N^\ast d^c(r^2)=\pi_N^\ast (2r^2d^c\log r)\stackrel{\eqref{contactOneForm}}{=}\pi_N^\ast (2r^2\tilde{\eta})=\beta$, see Section \ref{ConicalSect}.

Now we compute the one-forms $\theta_j^P$, $j=0,1,2,3$ on $P$, introduced
in \re{thetaEqu}:
\begin{align*}
\theta_0^P&=-\frac{1}{2}df=-rdr,\\
\theta_1^P&=\eta+\frac{1}{2}\beta=ds+\frac{1}{2}r^2\tilde{\eta}+\eta_{can}=ds+\frac{f+c}{2}\tilde{\eta}+\eta_{can},\\
\theta_2^P&=\frac{1}{2}\omega_3(Z,\cdot)=-\frac{i}{2}\sum(\zb^I \bar{A}_I-z^I A_I)=-\mathrm{Im} \sum z^IA_I,\\
\theta_3^P&=-\frac{1}{2}\omega_2(Z,\cdot)=\frac{1}{2}\sum (z^I A_I+\zb^I \bar{A}_I)=\mathrm{Re} \sum z^I A_I.\\
\end{align*}
For the calculation of $\theta_2^P$ and $\theta_3^P$, we used $Z=2i\sum(z^I\frac{\partial}{\partial z^I}-\zb^I\frac{\partial}{\partial \zb^I})^h$ and
\re{RigidCmapOmega2}-\re{RigidCmapOmega3}.

We compute the pseudo-Riemannian metric
\[g_P=\frac{2}{f_1}\eta^2+\pi^\ast g_N\stackrel{\eqref{RigidCmapMetric}}{=}\frac{2}{f_1}(ds+\frac{c}{2}\tilde{\eta}+\eta_{can}+\frac{f_1}{2}\tilde{\eta})^2+g_M+\sum A_IN^{IJ}\bar{A}_J\]
and the degenerate tensor field
\begin{align*}
\tilde{g}_P&=g_P-\frac{2}{f}\sum_{j=0}^{3}(\theta_j^P)^2\\&=g_P-\frac{2}{f}\left(r^2dr^2+(ds+\frac{c}{2}\tilde{\eta}+\eta_{can}+\frac{f}{2}\tilde{\eta})^2+(\sum z^I A_I)(\sum\zb^J\bar{A}_J)\right)\\
&=\left(\frac{2}{f_1}-\frac{2}{f}\right)(ds+\frac{c}{2}\tilde{\eta}+\eta_{can})^2+\left(\frac{f_1}{2}-\frac{f}{2}\right)\tilde{\eta}^2-\frac{2}{f}r^2dr^2+g_M\\&\quad\quad\quad\quad~
\,\quad\quad\quad\quad\quad\quad\quad\quad +\sum A_IN^{IJ}\bar{A}_J-\frac{2}{f}(\sum z^I A_I)(\sum\zb^J\bar{A}_J),
\end{align*}
see \re{gPEqu} and \re{gPtildeEqu}. 
As always, pullbacks from $M$ and $N$ to $P$ are implied where necessary.
Using $\frac{f_1}{2}-\frac{f}{2}=-r^2=-(f+c)$, $\frac{2}{f_1}-\frac{2}{f}=-\frac{4}{f}\frac{f+c}{f+2c}$, $\frac{2}{f}=\frac{2}{r^2}+\frac{2c}{f(f+c)}$ and $g_M\stackrel{\eqref{dKconeMetricDecomposition}}{=}dr^2+r^2(\tilde{\eta}^2-g_{\bar{M}})$,
we get
\begin{align*}
\tilde{g}_P&=-r^2g_{\bar{M}}-\frac{f+2c}{f}dr^2-\frac{4}{f}\frac{f+c}{f+2c}(ds+\frac{c}{2}\tilde{\eta}+\eta_{can})^2-\frac{2c}{f(f+c)}(\sum z^I A_I)(\sum\zb^J\bar{A}_J)\\&\quad\quad\quad\quad \quad\quad\quad\quad\quad\quad+\sum A_IN^{IJ}\bar{A}_J-\frac{2}{r^2}(\sum z^I A_I)(\sum\zb^J\bar{A}_J).
\end{align*}
We claim that the last two terms can be combined into $-\frac{1}{2}\sum dp_a\hat{H}^{ab}dp_b$, which appeared in the Ferrara-Sabharwal metric \eqref{FSmetric}.
This will be proven in the lemma below, see \re{relationFiberMetricsCmap}. 

We use the local coordinates
\[r=\sqrt{\sum z^IN_{IJ}\zb^J},~\phi:=\arg z^0,~X^\mu=\frac{z^\mu}{z^0}\]
on the conical affine special K\"ahler base $M$ and choose the submanifold $N'=\{\phi=0\}\subset P=N\times S^1$, which is transversal to 
\[Z_1^P=(Z-\eta(Z)X_P)+f_1 X_P=Z+(r^2+f_1)X_P=2\partial_\phi-c\partial_s,\]
where $X_P=\partial_s$ is the fundamental vector field on $P$, cf.\
\re{Z1PEqu}.

In these coordinates, we have
\[ |z^0|^2 = r^2 e^{\mathcal{K}}\] 
and, hence,  
\begin{equation*}
\tilde{\eta}=\frac12 d^c \log r^2 = \frac12 d^c \log |z^0|^2 -\frac{1}{2}d^c\mathcal{K}= d\phi-\frac{1}{2}d^c\mathcal{K}=d\phi+\sum\frac{iN_{IJ}(X)}{2X^tN\bar{X}}(X^Id\bar{X}^J-\bar{X}^JdX^I)
\end{equation*}
and
\begin{align*} \sum(z^IA_I)\sum(\zb^J \bar{A}_J)=|z^0|^2\sum (X^IA_I)\sum(\bar{X}^J\bar{A}_J)=r^2e^{\mathcal{K}}|\sum(X^Id\zt_I+F_I(X)d\z^I)|^2,
\end{align*}
where $\mathcal{K}=-\log X^tN\bar{X}$, $X^tN\bar{X} :=\sum X^IN_{IJ}\bar{X}^J$,  
is the K\"ahler potential for the projective special K\"ahler metric $g_{\bar{M}}$.
Replacing the coordinates $r$ and $s$ by $\rho:=f$ and $\tilde{\phi}:=-4s$ and
recalling that $\sigma =\mathrm{sgn}\,f$, we obtain the quaternionic K\"ahler metric $g'=\frac{1}{2|f|}\tilde{g}_P|_{N'}$ from the HK/QK correspondence (Theorem \ref{mainThm}) such that $g^c_{FS}:=-2\sigma g'$ is given by
\begin{align}g^c_{FS}=\frac{\rho+c}{\rho}g_{\bar{M}}&+\frac{1}{4\rho^2}\frac{\rho+2c}{\rho+c}d\rho^2+\frac{1}{4\rho^2}\frac{\rho+c}{\rho+2c}(d\tilde{\phi}+\sum(\z^Id\zt_I-\zt_Id\z^I)+cd^c\mathcal{K})^2\nonumber \\
&+\frac{1}{2\rho}\sum dp_a\hat{H}^{ab}dp_b+\frac{2c}{\rho^2}e^{\mathcal{K}}\left|\sum (X^Id\zt_I+F_I(X)d\z^I)\right|^2.\label{DefFSmetric}
\end{align}
For $c=0$, $g^c_{FS}$ reduces to the Ferrara-Sabharwal metric \eqref{FSmetric}.

Notice that the above metric $g^c_{FS}$ obtained from the HK/QK
correspondence is defined on a 
subset of $\bar{M}\times \bR^* \times S^1\times \bR^{2n+2}$, where 
the $\bR^*$-factor corresponds to the 
coordinate $\rho$ (which may now take negative values)
and the $S^1$-factor is parametrized by the coordinate
$\tilde{\phi}=-4s$ considered modulo 
$8\pi \bZ$. Replacing the above subset by
its universal covering (that is replacing 
$S^1$ by $\bR$) we obtain a subset of $\bar{M}\times \bR^* \times \bR^{2n+3}$.
In particular,  $g_{FS}=g^0_{FS}$ is defined
on $\bar{N}$ as well as
on the cyclic quotient $\bar{N}/\bZ = 
\bar{M}\times \bR^{>0} \times S^1\times \bR^{2n+2}$.   

The pseudo-hyper-K\"ahler metric $g_N$ has signature $(4,4n)$ and $Z$ is space-like. Hence, $g'$ is negative definite if $f>0$ and $f_1<0$, it has signature $(4,4n)$ if $f_1f>0$ and it has signature $(8,4(n-1))$ if $f<0$ and $f_1>0$ (see Corollary 1 in \cite{ACM}). Using $f=\rho$ and $f_1=-\rho-2c$, we get
\[\mathrm{sign}\,g'=\begin{cases} (0,4n+4) &\quad\quad\mbox{for } \rho>\max\{0,-2c\} \\
(4,4n) &\quad\quad\mbox{for } 0<\rho<-2c,~c<0
\\
(4,4n) &\quad\quad\mbox{for } -2c<\rho<0,~c>0
\\
(8,4(n-1)) &\quad\quad\mbox{for } \rho<\min\{0,-2c\}. \end{cases}\]
Taking into account that by definition $r^2=g_M(\xi,\xi)>0$, i.e. $\rho>-c$, we get
\[\mathrm{sign}\,g'=\begin{cases} (0,4n+4) &\quad\quad\mbox{for } \rho>\max\{0,-2c\}~(\Leftrightarrow r^2>|c|) \\
(4,4n) &\quad\quad\mbox{for } -c<\rho<\max\{0,-2c\}~(\Leftrightarrow 0<r^2<|c|). \end{cases}\]

It remains to prove
\bl 
\begin{equation}\label{relationFiberMetricsCmap}
\sum dp_a\hat{H}^{ab}dp_b=-2\sum A_{I}N^{IJ}\bar{A}_J+\frac{4}{r^2}(\sum z^IA_{I})(\sum \zb^J\bar{A}_J),
\end{equation}
where, as in the last section, $(p_a)=(\zt_I,\z^J)$ and $(\hat{H}^{ab})=\begin{pmatrix}\mathcal{I}^{-1} & \mathcal{I}^{-1}\mathcal{R} \\ \mathcal{R}\mathcal{I}^{-1} & \mathcal{I}+\mathcal{R}\mathcal{I}^{-1}\mathcal{R}\end{pmatrix}$.
\el
\pf 
Recall that $A_I=d\zt_I+\sum_J F_{IJ} d\z^J$, $I=0,\ldots,n$. We write $A=(A_I)=d\zt+\F d\z$, where $d\zt=(d\zt_I)$, $d\z=(d\z^I)$ are form-valued column vectors and $\F:=(F_{IJ})$.

First, we show that $\sum A_IN^{IJ}\bar{A}_J=\sum dp_aH^{ab}dp_b$ with
\[(H^{ab}):=\begin{pmatrix} N^{-1} & \frac{1}{2}N^{-1}R \\
  \frac{1}{2}RN^{-1} & \frac{1}{4}(N+RN^{-1}R)\end{pmatrix}, \]
where $R := 2\mathrm{Re}\, \F$: 
\begin{align*}
\sum A_IN^{IJ}\bar{A}_J&=(d\zt^t+d\z^t\F)N^{-1}(d\zt+\overline{\F}d\z)\\
&=(d\zt^t+d\z^t\frac{1}{2}(R+iN))N^{-1}(d\zt+\frac{1}{2}(R-iN)d\z)\\
&=d\zt^tN^{-1}d\z+d\zt^t\frac{1}{2}N^{-1}R\,d\z+d\z^t\frac{1}{2}RN^{-1}d\zt+d\z^t\frac{1}{4}(N+RN^{-1}R)d\z.
\end{align*}
Now, we show that $(\sum z^IA_I)(\sum \bar{A}_J\zb^J)=\sum dp_a \breve{H}^{ab}dp_b$ with
\[(\breve{H}^{ab}):=\frac{1}{2}\begin{pmatrix} z\zb^t+\zb z^t & z\zb^t\overline{\F}+\zb z^t\F \\ \overline{\F}\zb z^t+\F z\zb^t & \F z\zb^t\overline{\F}+\overline{\F}\zb z^t\F \end{pmatrix}: \]
\begin{align*}
(\sum z^I A_I)(\sum \zb^J\bar{A}_J)&=(d\zt^t z+d\z^t\, \F z)(\zb^td\zt + \zb^t\overline{\F} d\z)\\
&=d\zt^t z\zb^td\zt +d\zt^t z\zb^t\overline{\F} d\z+d\z^t \F z\zb^t d\zt+d\z^t\F z\zb^t \overline{\F} d\z\\
&=d\zt^t \frac{1}{2}(z\zb^t+\zb z^t)d\zt +d\zt^t \frac{1}{2}(z\zb^t\overline{\F}+\zb z^t\F) d\z\\
&\quad\quad\quad\quad+d\z^t \frac{1}{2}(\F z\zb^t+\overline{\F} \zb z^t) d\zt +d\z^t\frac{1}{2}(\F z\zb^t \overline{\F}+\overline{\F} \zb z^t \F) d\z.
\end{align*}
Hence, the right side of equation \eqref{relationFiberMetricsCmap} is given by $\sum dp_a (-2H^{ab}+\frac{4}{r^2}\breve{H}^{ab}) dp_b$.

To rewrite the left side of \eqref{relationFiberMetricsCmap}, we need to invert
$\mathcal{I}=\mathrm{Im}\,\mathcal{N}=-\frac{1}{2}N+\frac{Nzz^tN}{2z^tNz}+\frac{N\zb \zb^tN}{2\zb^tN\zb}$. It is easy to check that the inverse of $\mathcal{I}$ is given by \cite{MV}
\[\mathcal{I}^{-1}=-2N^{-1}+\frac{2}{z^tN\zb}(z\zb^t+\zb z^t).\]
Using $\mathcal{R}=\mathrm{Re}\,\mathcal{N}=\frac{1}{2}R+\frac{iNzz^tN}{2z^tNz}-\frac{iN\zb \zb^tN}{2\zb^tN\zb}$, we obtain
\begin{align*}
\mathcal{I}^{-1}\mathcal{R}=-N^{-1}R+\frac{1}{z^tN\zb}(z\zb^t(R-iN)+\zb z^t(R+iN))=-N^{-1}R+\frac{2}{r^2}(z\zb^t\overline{\F}+\zb z^t\F)
\end{align*}
and hence
\[ \mathcal{R}\mathcal{I}^{-1}=(\mathcal{I}^{-1}\mathcal{R})^t=-RN^{-1}+\frac{2}{r^2}(\overline{\F}\zb z^t+\F z \zb^t).\]
For the lower right block in $(\hat{H}^{ab})$, we calculate
\begin{align*}
\mathcal{R}\mathcal{I}^{-1}\mathcal{R}&=-\frac{1}{2}RN^{-1}R+\frac{1}{z^tN\zb}(\overline{\F}\zb z^t(R+iN)+\F z\zb^t (R-iN))\\
&\quad\quad\quad\quad\,\quad\quad +\frac{i}{z^tNz}(-\frac{1}{2}R+F)zz^t N-\frac{i}{\zb^tN\zb}(-\frac{1}{2}R+\overline{\F})\zb \zb^tN\\
&= -\frac{1}{2}RN^{-1}R +\frac{2}{r^2}(\overline{\F}\zb z^t\F+\F z \zb^t\overline{\F})-\frac{Nzz^tN}{2z^tNz}-\frac{N\zb \zb^t N}{2\zb^tN \zb}
\end{align*}
and hence
\[\mathcal{I}+\mathcal{R}\mathcal{I}^{-1}\mathcal{R}=-\frac{1}{2}(N+RN^{-1}R)+\frac{2}{r^2}(\overline{\F}\zb z^t\F+\F z \zb^t\overline{\F}).\]
This shows that $(\hat{H}^{ab})=-2(H^{ab})+\frac{4}{r^2}(\breve{H}^{ab})$ and thus proves Eq.\ \eqref{relationFiberMetricsCmap}.
\qed

This proves Theorem \ref{HKQKThm}.\qed  

\br
\rm{Note that the quaternionic K\"ahler metric $g_{FS}^c$ given in \eqref{DefFSmetric} agrees with the one-loop deformed Ferrara-Sabhwarwal metric first obtained in \cite{RSV} (see also \cite{APP}, Eq. (2.93)).}
\er

\br \label{RemarkKK}
\rm{One can check that the restriction of the metric $g'=-\frac{\sigma}{2}
g^{c}_{FS}$ given in \eqref{DefFSmetric} to $M':=\{\z=\zt=0\}\subset
N'$ is the metric that one obtains when applying the K/K
correspondence to the original conical affine special K\"ahler
manifold $(M,J,g_M,\nabla,\xi)$ with respect to the holomorphic
Killing field $Z=2J\xi$. This is a special case of 
Theorem \ref{thmKKexample}.
For $c=0$, we have \[(M',g_{FS}|_{M'})=(\mathbb{R}^{>0}\times S^1\times \bar{M},g_{\mathbb{C}H^1}+g_{\bar{M}}),\] where $(\bar{M},J_{\bar{M}},g_{\bar{M}})$ is the underlying projective special K\"ahler manifold and $(\mathbb{R}^{>0}\times S^1,g_{\mathbb{C}H^1})$ is a $\mathbb{Z}$-quotient of the complex hyperbolic line. The hyperbolic metric $g_{\mathbb{C}H^1}=\frac{1}{4\rho^2}(d\rho^2+d\tilde{\phi}^2)$ is normalized such that its scalar curvature is $-8$.} 
\er

\begin{appendix} \section{A simple example of the HK/QK 
correspondence}
Here we consider 
\[ M=\{ (z,w)\in \bC^2| z\neq 0\} \subset \bC^2\] 
with its standard
flat hyper-K\"ahler 
structure\footnote{The flat hyper-K\"ahler manifold 
$\bC^2$ is also considered in \cite{Hi}
but the corresponding quaternionic K\"ahler metric
is not computed there.} $(g,J_1,J_2,J_3)$. 
This manifold 
is in the image of the
rigid c-map and therefore admits a vector 
field $Z$ verifying the above assumptions, see 
\cite{ACM}, Proposition 2. The canonical choice 
of function $f$ which leads to a definite
quaternionic K\"ahler metric $g'$ is $f=\frac{1}{4}g(Z,Z)$, 
see \cite{ACM}, Corollary 4. 
The metric $g'$ is in fact negative definite if we take $g$ positive
definite and vice versa. 

Let us first compute all the relevant geometric data on $M$ in terms 
of the standard $J_1$-holomorphic coordinates $(z,w)$ of $\bC^2$, which satisfy 
$J_3^*dw=d\bar{z}$. The metric and K\"ahler forms are given by: 
\begin{eqnarray*}
g&=& 2(|dz|^2 + |dw|^2),\\
\o_1&=&i(dz\wedge d\bar{z} + dw\wedge d\bar{w}),\\
\o_2 &=& i(dz\wedge dw - d\bar{z}\wedge d\bar{w}),\\
\o_3 &=& dz\wedge dw + d\bar{z}\wedge d\bar{w}. 
\end{eqnarray*}
The vector field $Z$ is given by
\begin{eqnarray*} Z &=& 2(iz\frac{\partial}{\partial z} - 
i \bar{z}\frac{\partial}{\partial 
\bar{z}}),\\
g(Z,Z) &=& 8|z|^2, 
\end{eqnarray*}
with the canonical choice of Hamiltonian given by 
\[ f= 2|z|^2,\] 
such that $df = -\o_1 Z$ and
\[ f_1 = f-\frac{g(Z,Z)}{2} =-2|z|^2.\]
Notice that the functions $f$ and $f_1$ are nowhere
vanishing on $M$.  
Then we consider the trivial $S^1$-principal bundle 
\[ P=M\times S^1,\quad  S^1=\{ e^{is}|s\in \bR\}.\] 
with the connection form  
\[ \eta = ds + \eta_M,\]
where $s$ is the natural coordinate on 
$S^1=\{ e^{is}|s\in \bR\}$ 
and $\eta_M$ is the following one-form on $M$ 
\[ \eta_M = \frac{i}{2}(zd\bar{z}-\bar{z}dz+wd\bar{w}-\bar{w}dw) -
\frac{1}{2}gZ.\]
Computing 
\be \label{gzEqu} gZ= 2i(zd\bar{z}-\bar{z}dz),\ee
we get 
\[ \eta_M = \frac{i}{2}(-zd\bar{z}+\bar{z}dz+wd\bar{w}-\bar{w}dw).\]

\noindent 
\br\rm{Notice that, in the trivialization of $P$, $X_P=\partial_s$, 
$\tilde{Z}= Z -\eta_M(Z)\partial_s$ and $Z_1 = Z +(-\eta_M(Z)+f_1)
\partial_s$. The above formula for $\eta_M$ implies that 
$\eta_M(Z)=-2|z|^2=f_1$ and, thus, $Z_1=Z$, 
in the given trivialization.}
\er

We define $M'$ as the submanifold of $P=M\times S^1$ defined by
$\mathrm{Im}\,z=0$. We will use $w$, $r = \sqrt{2}|z|$ and $s$ as local
coordinates on $M'$. $M'$ intersects each orbit of the
$S^1$-action $S^1_{Z_1}$ generated by $Z_1$ in exactly one point  
such that we can identify $M'$ with the orbit space
$P/S^1_{Z_1}$. Now we compute the one-forms $gZ$ and 
$\theta_a$, $a=0,1,2,3$, on $M'$. 
{}Writing $z=\frac{r}{\sqrt{2}} e^{i\arg z}$, from \re{gzEqu} we get 
\[ gZ = 2r^2 d\arg z ,\]
which shows that $gZ$ vanishes on $M'$ and that 
\[ \eta_M|_{M'} =  \frac{i}{2}(wd\bar{w}-\bar{w}dw) = \frac{1}{4}(\zt d\z-\z d\zt)
= \eta_{can},\]
if we write $w=\frac{1}{2}(\zt+i\z)$ and define $\eta_{can}=\frac{1}{4}(\zt d\z-\z d\zt)$. Using 
that $gZ|_{M'}=0$, we have: 
\begin{eqnarray*}
\theta_0|_{M'} &=& -(zd\bar{z}+\bar{z}dz)|_{M'}=
-r dr,\\
\theta_1|_{M'} &=& \eta|_{M'},\\
\theta_2|_{M'}  &=& \frac{i}{\sqrt{2}}r (dw-d\bar{w}),\\
\theta_3|_{M'}  &=&\frac{1}{\sqrt{2}}r (dw+d\bar{w}),
\end{eqnarray*}
which implies 
\[ \sum_a (\theta_a^P)^2 = (\eta|_{M'})^2+ r^2 (dr^2 + 2|dw|^2).\] 
So 
\[ \tilde{g}_P|_{M'} = 
\left( \frac{2}{f_1}-\frac{2}{f}\right) \eta^2|_{M'} + (dr^2 + 2|dw|^2)
-2(dr^2 + 2|dw|^2)= \left( \frac{2}{f_1}-\frac{2}{f}\right)\eta^2|_{M'}
-(dr^2 + 2|dw|^2).\] 
Now 
\[ \frac{2}{f_1}-\frac{2}{f}= \frac{2(f-f_1)}{ff_1}=
\frac{g(Z,Z)}{ff_1}= -\frac{4}{r^2}\]
and 
\[\eta|_{M'} = ds + \eta_{can}.\]
Therefore, we can rewrite
\[ \tilde{g}_P|_{M'} = -\frac{4}{r^2}(ds+\eta_{can})^2-(dr^2 + 
2|dw|^2)\] 
and 
\[ -2g' = -\frac{1}{f}\tilde{g}_P|_{M'} = \frac{4}{r^4}
( ds+\eta_{can})^2
+\frac{dr^2}{r^2} + 
2\frac{|dw|^2}{r^2}.\]
Putting $\rho= r^2$ and $\tilde{\phi}=-4s$, we can rewrite this as
\[ -2g'= \frac{1}{4\rho^2}(d\tilde{\phi}+\z d\zt-\zt d\z)^2 +\frac{d\rho^2}{4\rho^2} + 
\frac{d\zt^2+d\z^2}{2\rho} .\] 
The Riemannian metric $-2g'$ is precisely the 
Ferrara-Sabharwal metric (cf. \cite{FS,CHM}),
which in the present case coincides with the complex hyperbolic
metric. 

\noindent
\br \rm{The manifold $M'$ is a cyclic quotient of 
the complex hyperbolic plane. The complex hyperbolic plane 
is  parametrized by the global coordinates $(w,r >0, s)$ 
and $M'$ is obtained as the quotient by $2\pi \bZ$ acting 
by translations in $s$.} 
\er

\br \label{Remark4D}
\rm{Carrying out the above calculation with the Hamiltonian function $f$ replaced by $f-c$, one obtains the deformed Ferrara-Sabharwal metric \eqref{DefFSmetric} for the special case where the underlying projective special K\"ahler manifold is a point (i.e. for holomorphic prepotential $F=\frac{i}{2}(z^0)^2$):
\[g_{UH}^c:=\frac{1}{4\rho^2}\left( \frac{\rho+2c}{\rho+c}d\rho^2+\frac{\rho+c}{\rho+
2c}(d\tilde{\phi}+\zeta^0 d\tilde{\zeta}_0-\tilde{\zeta}_0d\z^0)^2+2(\rho+2c)((d\zt_0)^2+(d\z^0)^2)\right).\]
This is known to physicists as the one-loop corrected universal hypermultiplet metric and was derived in \cite{AMTV}. As was already noticed in \cite{AMTV}, this metric admits an isometric action of the three-dimensional Heisenberg group generated by the Killing vector fields
\[\frac{\partial}{\partial \tilde{\phi}},\quad \frac{\partial}{\partial \tilde{\zeta_0}}+\zeta^0\frac{\partial}{\partial \tilde{\phi}},\quad \frac{\partial}{\partial \zeta^0}-\tilde{\zeta}_0\frac{\partial}{\partial \tilde{\phi}};\]
and hence falls under the classification of 4-dimensional self-dual
Einstein metrics with non-zero scalar curvature admitting two
commuting Killing vector fields by Calderbank and
Pedersen\footnote{Calderbank and Pedersen express such a metric in
  terms of an eigenfunction $\mathcal{F}(r,\eta)$ of the Laplacian on
  the hyperbolic plane $\{(r,\eta)\in
  \mathbb{R}^{>0}\times\mathbb{R}\}$. In their formalism, the metric
  $-2g_{UH}^c$ corresponds to the function $\mathcal{F}=\frac{r^2-c}{\sqrt{r}}$ and their coordinates $(r,\eta,\phi,\psi)$ are related to our coordinates by $\rho=r^2-c,~\tilde{\phi}=2\psi+\phi\eta,~\zt_0=\frac{1}{\sqrt{2}}\phi,~\z^0=\sqrt{2}\eta$.} \cite{CP}.

\noindent For $c>0$, $g^c_{UH}$ is positive definite and of negative
scalar curvature on the domains $\{-c<\rho<0\}$ and $\{\rho>0\}$ in
$\mathbb{R}^4$. For $c<0$, $g^c_{UH}$ is positive definite and of
negative scalar curvature on $\{\rho>-2c\}\subset \mathbb{R}^4$ and
$-g^c_{UH}$ is positive definite and of positive scalar curvature on
$\{-c<\rho<-2c\}\subset \mathbb{R}^4$, cf.\ Theorem \ref{HKQKThm}.

\noindent Notice that the complex hyperbolic metric $g_{UH}^0$ is
symmetric and hence complete. Using this, one can show that $g^c_{UH}$
is complete on the domain $\{\rho>0\}$ for $c>0$. In fact,
$g^c_{UH}> \frac12 g_{UH}^0$. 
On the other domains
mentioned above, however, the positive or negative definite metric
$g^c_{UH}$ is incomplete, as stated in the proposition below. 
This is in agreement with the result of A.\ Haydys, who
studied $-g^c_{UH}$ on $\{-c<\rho<-2c\}\subset \mathbb{R}^4$ for the
special case $c=-1$ (see \cite{Haydys}, Example 9 (resp. 3.2 in the
arXiv version)).}
\er 
\bp 
\begin{enumerate}
\item[(i)]  For $c>0$  the quaternionic K\"ahler metric $g^c_{UH}$ of 
negative scalar curvature is complete 
on the domain $\{\rho>0\}$ and incomplete on $\{-c<\rho<0\}$. 
\item[(ii)] For $c<0$  the quaternionic K\"ahler metric $g^c_{UH}$ is 
incomplete and of negative scalar curvature on $\{\rho>-2c\}$. 
\item[(iii)] For $c<0$  the quaternionic K\"ahler metric $-g^c_{UH}$ is 
incomplete and of positive scalar curvature on $\{-c<\rho<-2c\}$. 
\end{enumerate}
\ep 
\pf 
It remains to prove the incompleteness in the corresponding cases. 
In cases (i) and (iii) we consider the 
curve 
\[ \r = t-c,\quad \tilde{\phi}=\tilde{\zeta}_0=\zeta^0=0,\quad
0<t<\frac{|c|}{2},\]
which approaches the boundary of the respective domain 
for $t\ra 0$. Its length is
given by 
\[ \frac12 \int_0^{\frac{|c|}{2}}\frac{1}{|t-c|}\sqrt{\frac{|t+c|}{t}}dt
\ge C\int_0^{\frac{|c|}{2}}\frac{dt}{2\sqrt{t}}< \infty,
 \]
where $C>0$ is a lower bound for the 
continuous function $\frac{\sqrt{|t+c|}}{|t-c|}$ on the compact
interval $[0,\frac{|c|}{2}]$. 
In case (ii) we consider instead the curve 
\[ \r = t-2c,\quad \tilde{\phi}=\tilde{\zeta}_0=\zeta^0=0,\quad
0<t<1,\]
which approaches the boundary for $t\ra 0$. 
Its length is the integral of the continuous 
function $\frac{1}{2|t-2c|}\sqrt{\frac{t}{t-c}}$
on the compact
interval $[0,1]$ and, hence, finite. 
\qed
\br \rm{Notice that the above
proof for the incompleteness in case (ii) is still valid
in higher dimensions for the positive definite quaternionic K\"ahler
metric $g^c_{FS}$, $c<0$, 
on the domain $\{ \r > -2c\}$, see Theorem \ref{HKQKThm} for a 
description of the domain of positivity of the one-loop
deformed Ferrara-Sabharwal metric \re{DefFSmetric} depending on the sign of
$c$. On the contrary, the proof of the completeness in case
(i), given in Remark \ref{Remark4D}, does not extend in a straightforward way to higher
dimensions.\footnote{Note added in proof: In the meantime it has been proven that the one-loop deformed 
Ferrara-Sabharwal metric $g^c_{FS}$ for $c>0$ on the domain $\rho >0$ is  complete for every complete 
projective special K\"ahler manifold in the image of the supergravity $r$-map, see M.\ Dyckmanns, PhD dissertation, University of Hamburg, to appear in 2015.}}     
\er  
\end{appendix}

\end{document}